\documentclass[draft,12pt,amssymb,amstex]{article}

\usepackage{amsmath}
\usepackage{amssymb}
\usepackage{amsxtra}
\usepackage{latexsym}
\usepackage{ifthen}

\usepackage{mathtext}
\usepackage[cp1251]{inputenc}
\usepackage[T2A]{fontenc}
\usepackage[dvips]{graphicx}

\begin{document}

\newcounter{lemma}[section]
\newcommand{\lemma}{\par \refstepcounter{lemma}%
{\bf Lemma \arabic{section}.\arabic{lemma}.}}
\renewcommand{\thelemma}{\thesection.\arabic{lemma}}

\newcounter{corol}[section]
\newcommand{\corol}{\par \refstepcounter{corol}%
{\bf Corollary \arabic{section}.\arabic{corol}.}}
\renewcommand{\thecorol}{\thesection.\arabic{corol}}

\newcounter{rem}[section]
\newcommand{\rem}{\par \refstepcounter{rem}%
{\bf Remark \arabic{section}.\arabic{rem}.}}
\renewcommand{\therem}{\thesection.\arabic{rem}}

\newcounter{theo}[section]
\newcommand{\theo}{\par \refstepcounter{theo}%
{\bf Theorem \arabic{section}.\arabic{theo}.}}
\renewcommand{\thetheo}{\thesection.\arabic{theo}}

\newcounter{propo}[section]
\newcommand{\propo}{\par \refstepcounter{propo}%
{\bf Proposition \arabic{section}.\arabic{propo}.}}
\renewcommand{\thepropo}{\thesection.\arabic{propo}}

\numberwithin{equation}{section}

\newcommand{\osc}{\operatornamewithlimits{osc}}

\def\Xint#1{\mathchoice
   {\XXint\displaystyle\textstyle{#1}}%
   {\XXint\textstyle\scriptstyle{#1}}%
   {\XXint\scriptstyle\scriptscriptstyle{#1}}%
   {\XXint\scriptscriptstyle\scriptscriptstyle{#1}}%
   \!\int}
\def\XXint#1#2#3{{\setbox0=\hbox{$#1{#2#3}{\int}$}
     \vcenter{\hbox{$#2#3$}}\kern-.5\wd0}}
\def\dashint{\Xint-}

\markboth{\centerline{ R. SALIMOV
 }} {\centerline{LOWER $Q$-HOMEOMORPHISMS WITH RESPECT TO $p$-MODULUS AND ORLICZ-SOBOLEV CLASSES}}

\def\cc{\setcounter{equation}{0}
\setcounter{figure}{0}\setcounter{table}{0}}

\overfullrule=0pt

\title{{\bf LOWER $Q$-HOMEOMORPHISMS WITH RESPECT TO $p$-MODULUS AND ORLICZ-SOBOLEV CLASSES}}

\author{\bf  R. Salimov}

\maketitle

\Large \abstract
We show that under a condition of the Calderon type on
$\varphi$  the homeomorphisms $f$ with finite distortion in
$W^{1,\varphi}_{\rm loc}$ and, in particular, $f\in W^{1,s}_{\rm
loc}$ for $s>n-1$ are the so-called lower $Q$-homeomorphisms  with respect to $p$-modulus where
$Q(x)$ is equal to its outer $p$-dilatation $K_{p,f}(x)$.  \endabstract

\bigskip



\large

\section{Introduction}

In what follows, $D$ is a domain in a finite-dimensional Euclidean
space. Following Orlicz, see \cite{Or1}, given a convex increasing
function $\varphi:[0,\infty)$ $\to[0,\infty)$, $\varphi(0)=0$,
denote by $L_{\varphi}$ the space of all functions $f:D\to{\Bbb
R}$ such that
\begin{equation}\label{eqOS1.1}\int\limits_{D}\varphi\left(\frac{|f(x)|}
{\lambda}\right)\,dm(x)<\infty\end{equation} for some $\lambda>0$
where $dm(x)$ corresponds to the Lebesgue measure in $D$.
$L_{\varphi}$ is called the {\bf Orlicz space}.  If
$\varphi(t)=t^p$, then we write also $L_{p}$. In other words,
$L_{\varphi}$ is the cone over the class of all functions
$g:D\to{\Bbb R}$ such that
\begin{equation}\label{eqOS1.2}
\int\limits_{D}\varphi\left(|g(x)|\right)\,dm(x)<\infty\end{equation}
which is also called the {\bf Orlicz class}, see \cite{BO}.

\medskip

The {\bf Orlicz-Sobolev class} $W^{1,\varphi}_{\rm loc}(D)$ is the
class of locally integrable functions $f$ given in $D$ with the
first distributional derivatives whose gradient $\nabla f$ belongs
locally in $D$ to the Orlicz class. Note that by definition
$W^{1,\varphi}_{\rm loc}\subseteq W^{1,1}_{\rm loc}$. As usual, we
write $f\in W^{1,p}_{\rm loc}$ if $\varphi(t)=t^p$, $p\geqslant1$.
It is known that a continuous function $f$ belongs to
$W^{1,p}_{\rm loc}$ if and only if $f\in ACL^{p}$, i.e., if $f$ is
locally absolutely continuous on a.e. straight line which is
parallel to a coordinate axis, and if the first partial
derivatives of $f$ are locally integrable with the power $p$, see,
e.g., 1.1.3 in \cite{Maz}. The concept of the distributional
derivative was introduced by Sobolev \cite{So} in ${\Bbb R}^n$,
$n\geqslant2$, and it is developed under wider settings at
present, see, e.g., \cite{Re1}.

\medskip

Later on, we also write $f\in W^{1,\varphi}_{\rm loc}$ for a locally
integrable vector-function $f=(f_1,\ldots,f_m)$ of $n$ real
variables $x_1,\ldots,x_n$ if $f_i\in W^{1,1}_{\rm loc}$ and
\begin{equation}\label{eqOS1.2a}
\int\limits_{D}\varphi\left(|\nabla
f(x)|\right)\,dm(x)<\infty\end{equation} where $|\nabla
f(x)|=\sqrt{\sum\limits_{i,j}\left(\frac{\partial f_i}{\partial
x_j}\right)^2}$. Note that in this paper we use the notation
$W^{1,\varphi}_{\rm loc}$ for more general functions $\varphi$ than
in the classical Orlicz classes giving up the condition on convexity
of $\varphi$. Note also that the Orlicz--Sobolev classes are
intensively studied in various aspects at present.

\medskip

Recall that a homeomorphism $f$ between domains $D$ and $D'$ in
${\Bbb R}^n$, $n\geqslant2$, is called of {\bf finite distortion}
if $f\in W^{1,1}_{\rm loc}$ and
\begin{equation}\label{eqOS1.3}\Vert
f'(x)\Vert^n\leqslant K(x)\cdot J_f(x)\end{equation} with a.e.
finite function $K$ where $\Vert f'(x)\Vert$ denotes the matrix norm
of the Jacobian matrix $f'$ of $f$ at $x\in D$,
$||f'(x)||=\sup\limits_{h\in{\Bbb R}^n,|h|=1}|f'(x)\cdot h|$, and
$J_f(x)={\rm det}f'(x)$ is its Jacobian.  We set $K_{p,f}(x)=\Vert f'(x)\Vert^p/$ $J_f(x)$
if $J_f(x)\ne 0,$ $K_{p,f}(x)=1$ if $f'(x)=0$ and $K_{p,f}(x)=\infty$ at the
rest points.

\medskip

First this notion was introduced on the plane for $f\in
W^{1,2}_{\rm loc}$ in the work \cite{IS}. Later on, this condition
was changed by $f\in W^{1,1}_{\rm loc}$ but with the additional
condition $J_f\in L^1_{\rm loc}$ in the monograph \cite{IM}. The
theory of the mappings with finite distortion had many successors,
see, e.g., a number of references in the monograph \cite{MRSY}.
They had as predecessors the mappings with bounded distortion, see
\cite{Re} and \cite{Vo}, in other words, the quasiregular
mappings, see, e.g., \cite{BI}, \cite{BH}, \cite{HKM}, \cite{MRV},
\cite{Ri} and \cite{Vu}.

\medskip

Note that the above additional condition $J_f\in L^1_{\rm loc}$ in
the definition of the mappings with finite distortion can be omitted
for ho\-meo\-mor\-phisms. Indeed, for each homeomorphism $f$ between
domains $D$ and $D'$ in ${\Bbb R}^n$ with the first partial
derivatives a.e. in $D$, there is a set $E$ of the Lebesgue measure
zero such that $f$ satisfies $(N)$-property by Lusin on $D\setminus
E$ and
\begin{equation}\label{eqOS1.1.1}\int\limits_{A}J_f(x)\,dm(x)=|f(A)|\end{equation}
for every Borel set $A\subset D\setminus E$, see, e.g., 3.1.4,
3.1.8 and 3.2.5 in \cite{Fe}. On this base, it is easy by the
H\"older inequality to verify, in particular, that if $f\in
W^{1,1}_{\rm loc}$ is a ho\-meo\-mor\-phism and $K_f\in L^q_{\rm
loc}$ for some $q>n-1$, then also $f\in W^{1,p}_{\rm loc}$ for
some $p>n-1$, that we use further to obtain corollaries.

\medskip

In this paper $H^k(A)$, $k\geqslant 0$, ${\rm dim}_H A$ denote the
{\bf $\bf k$-dimensional Hausdorff measure} and the {\bf Hausdorff
dimension}, correspondingly, of a set $A$ in ${\Bbb R}^n$,
$n\geqslant1$. It was shown in \cite{GV} that a set $A$ with
$\mbox{dim}_H\, A =p$ can be transformed into a set $B=f(A)$ with
$\mbox{dim}_H\, B =q$ for each pair of numbers $p$ and $q\in
(0,n)$ under a quasiconformal mapping $f$ of ${\Bbb R}^n$ onto
itself, cf. also \cite{Ba} and \cite{Bi$_2$}.

\cc
\section{Preliminaries}

First of all, the following fine property of functions $f$ in the
Sobolev classes $W^{1,p}_{\rm loc}$ was proved in the monograph
\cite{GR}, Theorem 5.5, and can be extended to the Orlicz-Sobolev
classes. The statement follows directly from the Fubini theorem
and the known characterization of functions in Sobolev's class
$W^{1,1}_{\rm loc}$ in terms of ACL (absolute continuity on
lines), see, e.g., Section 1.1.3 in \cite{Maz}.

\medskip

\medskip

\medskip

\medskip

\begin{theo}{}\label{thOS4.0} {\it Let $\Omega$ be an open set in ${\Bbb R}^{n}$,
$n\geqslant3$, and let $f:\Omega\to{\Bbb R}^{n}$ be a continuous
open mapping in the class $W^{1,\varphi}_{\rm loc}(\Omega)$ where
$\varphi:[0,\infty)\to[0,\infty)$ is increasing with the condition
\begin{equation}\label{eqOS3.0a}\int\limits_{1}^{\infty}\left[\frac{t}{\varphi(t)}\right]^
{\frac{1}{n-2}}dt<\infty.\end{equation}
. Then $f$ has a total differential a.e. in
$\Omega$.} \end{theo}

\medskip

\begin{corol}{}\label{corOS4.2b} {\it If $f:\Omega\to{\Bbb R}^{n}$ is a homeomorphism in
$W^{1,1}_{\rm loc}$ with $K_f\in L^p_{\rm loc}$ for $p>n-1$, then
$f$ is differentiable a.e.} \end{corol}

\cc

\begin{theo}{}\label{thOS3.1} {\it Let $U$ be an open set in ${\Bbb R}^n$,
$n\geqslant3$, and let $\varphi:[0,\infty)\to[0,\infty)$ is
increasing with the condition (\ref{eqOS3.0a}). Then each
continuous mapping $f:U\to{\Bbb R}^m$, $m\ge 1$, in the class
$W^{1,\varphi}_{\rm loc}$ has the $(N)$-property (furthermore, it
is locally absolutely continuous) with respect to the
$(n-1)$-dimensional Hausdorff measure on a.e. hyperplane
$\mathcal{P}$ which is parallel to a fixed coordinate hyperplane
${\mathcal P}_0$. Moreover, $H^{n-1}(f(E))=0$ whenever $|\nabla
f|=0$ on $E\subset\mathcal{P}$ for a.e. such $\mathcal{P}$.}
\end{theo}

\medskip

Note that, if the condition (\ref{eqOS3.0a}) holds for an
increasing function $\varphi$, then the function
$\varphi_*=\varphi(c\,t)$ for $c>0$ also satisfies
(\ref{eqOS3.0a}). Moreover, the Hausdorff measures are
quasi-invariant under quasi-isometries. By the Lindel\"of property
of ${\Bbb R}^n$, $U\setminus\{x_0\}$ can be covered by a countable
collection of open segments of spherical rings in
$U\setminus\{x_0\}$ centered at $x_0$ and each such segment can be
mapped onto a rectangular oriented segment of ${\Bbb R}^n$ by some
quasi-isometry, see, e.g., I.5.XI in \cite{Ku} for the Lindel\"of
theorem. Thus, applying piecewise Theorem \ref{thOS3.1}, we obtain
the following.

\medskip

\begin{corol}\label{corOS3.2} {\it Under (\ref{eqOS3.0a}) each
$f\in W^{1,\varphi}_{\rm loc}$ has the $(N)$-property (furthermore,
it is locally absolutely continuous) on a.e. sphere $S$ centered at
a prescribed point $x_0\in{\Bbb R}^n$. Moreover, $H^{n-1}(f(E))=0$
whenever $|\nabla f|=0$ on $E\subseteq S$ for a.e. such sphere $S$.}
\end{corol}

\section{Moduli of families of surfaces}\label{12}

The recent development of the moduli method in the connection with
modern classes of mappings can be found in the monograph \cite{MRSY}
and further references therein.

Let $\omega$ be an open set in $\overline{{\Bbb R}^k}$,
$k=1,\ldots,n-1$. A (continuous) mapping $S:\omega\to{\Bbb R}^n$ is
called a $k$-dimensional surface $S$ in ${\Bbb R}^n$. Sometimes we
call the image $S(\omega)\subseteq{\Bbb R}^n$ the surface $S$, too.
The number of preimages \begin{equation}\label{eq8.2.3} N(S,y)={\rm
card}\,S^{-1}(y)={\rm card}\,\{x\in\omega:S(x)=y\},\ y\in{\Bbb
R}^n\end{equation} is said to be a {\bf multiplicity function} of
the surface $S$. In other words, $N(S,y)$ denotes the multiplicity
of covering of the point $y$ by the surface $S$. It is known that
the multiplicity function is lower semicontinuous, i.e.,
$$N(S,y)\geqslant\liminf_{m\to\infty}\:N(S,y_m)$$ for every sequence
$y_m\in{\Bbb R}^n$, $m=1,2,\ldots\,$, such that $y_m\to y\in{\Bbb
R}^n$ as $m\to\infty$, see e.g. \cite{RR$^*$}, p. 160. Thus, the
function $N(S,y)$ is Borel measurable and hence measurable with
respect to every Hausdorff measure $H^k$; see e.g. \cite{Sa}, p. 52.

Recall that a $k$-dimensional Hausdorff area in ${\Bbb R}^n$ (or
simply {\bf area}) associated with a surface $S:\omega\to{\Bbb R}^n$
is given by \begin{equation}\label{eq8.2.4}
{\cal{A}}_S(B)={\cal{A}}^{k}_S(B):=\int\limits_B N(S,y)\,dH^{k}y
\end{equation} for every Borel set $B\subseteq{\Bbb R}^n$ and,
more generally, for an arbitrary set that is measurable with respect
to $H^k$ in ${\Bbb R}^n$, cf. 3.2.1 in \cite{Fe}. The surface $S$ is
called {\bf rectifiable} if ${\cal{A}}_S({\Bbb R}^n)<\infty$, see
9.2 in \cite{MRSY}.

If $\varrho:{\Bbb R}^n\to[0,\infty]$ is a Borel function, then its
{\bf integral over} $S$ is defined by the equality
\begin{equation}\label{eq8.2.5} \int\limits_S \varrho\,d{\cal{A}}:=
\int\limits_{{\Bbb R}^n}\varrho(y)\,N(S,y)\,dH^ky\,.\end{equation}
Given a family $\Gamma$ of $k$-dimensional surfaces $S$, a Borel
function $\varrho:{\Bbb R}^n\to[0,\infty]$ is called {\bf
admissible} for $\Gamma$, abbr. $\varrho\in{\rm adm}\,\Gamma$, if
\begin{equation}\label{eq8.2.6}\int\limits_S\varrho^k\,d{\cal{A}}\geqslant1\end{equation}
for every $S\in\Gamma$. Given $p\in(0,\infty)$, the {\bf
$p$-modulus} of $\Gamma$ is the quantity
\begin{equation}\label{eq8.2.7} M_p(\Gamma)=\inf_{\varrho\in{\rm adm}\,\Gamma}
\int\limits_{{\Bbb R}^n}\varrho^p(x)\,dm(x)\,.\end{equation} We also
set \begin{equation}\label{eq8.2.8}
M(\Gamma)=M_n(\Gamma)\end{equation} and call the quantity
$M(\Gamma)$ the {\bf modulus of the family} $\Gamma$. The modulus is
itself an outer measure in the space of all $k$-dimensional
surfaces.

\bigskip

We say that $\Gamma_2$ is {\bf minorized} by $\Gamma_1$ and write
$\Gamma_2>\Gamma_1$ if every $S\subset\Gamma_2$ has a subsurface
that belongs to $\Gamma_1$. It is known that $M_p(\Gamma_1)\geqslant
M_p(\Gamma_2)$, see \cite{Fu}, p.~176-178. We also say that a
property $P$ holds for {\bf $p$-a.e.} (almost every) $k$-dimensional
surface $S$ in a family $\Gamma$ if a subfamily of all surfaces of
$\Gamma$, for which $P$ fails, has the $p$-modulus zero. If $0<q<p$,
then $P$ also holds for $q$-a.e. $S$, see Theorem 3 in \cite{Fu}. In
the case $p=n$, we write simply a.e.

\bigskip

\begin{rem}\label{rmk8.2.9} The definition of the modulus immediately
implies that, for every $p\in(0,\infty)$ and $k=1,\ldots,n-1$

\medskip

\noindent (1) $p$-a.e. $k$-dimensional surface in ${\Bbb R}^n$ is
rectifiable,

\medskip

\noindent (2) given a Borel set $B$ in ${\Bbb R}^n$ of (Lebesgue)
measure zero, \begin{equation}\label{eq8.2.10} {\cal
A}_S(B)=0\end{equation} for $p$-a.e. $k$-dimensional surface $S$ in
${\Bbb R}^n$. \end{rem}

\bigskip

The following lemma was first proved in \cite{KR$_1$}, see also
Lemma 9.1 in \cite{MRSY}.

\bigskip

\begin{lemma}\label{lem8.2.11} {\it Let $k=1,\ldots,n-1$, $p\in[k,\infty)$,
and let $C$ be an open cube in ${\Bbb R}^n$, $n\geqslant2$, whose
edges are parallel to coordinate axis. If a property $P$ holds for
$p$-a.e. $k$-dimensional surface $S$ in $C$, then $P$ also holds for
a.e. $k$-dimensional plane in $C$ that is parallel to a
$k$-dimensional coordinate plane $H$.} \end{lemma}

\bigskip

The latter a.e. is related to the Lebesgue measure in the
correspon\-ding $(n-k)$-dimensional coordinate plane $H^{\perp}$
that is perpendicular to $H$.

\bigskip

The following statement, see Theorem 2.11 in \cite{KR$_2$} or
Theorem 9.1 in \cite{MRSY}, is an analog of the Fubini theorem, cf.
e.g. \cite{Sa}, p. 77. It extends Theorem 33.1 in \cite{Va}, cf.
also Theorem 3 in \cite{Fu}, Lemma 2.13 in \cite{MRSY$_4$} and Lemma
8.1 in \cite{MRSY}.

\bigskip

\begin{theo}\label{th8.2.12} {\it Let $k=1,\ldots,n-1$, $p\in[k,\infty)$,
and let $E$ be a subset in an open set $\Omega\subset{\Bbb R}^n$,
$n\geqslant2$. Then $E$ is measurable by Lebesgue in ${\Bbb R}^n$ if
and only if $E$ is measurable with respect to area on $p$-a.e.
$k$-dimensional surface $S$ in $\Omega$. Moreover, $|E|=0$ if and
only if \begin{equation}\label{eq8.2.13}{\cal
A}_S(E)=0\end{equation} on $p$-a.e. $k$-dimensional surface $S$ in
$\Omega$.} \end{theo}

\bigskip

\begin{rem}\label{rmk8.2.14} Say by the Lusin theorem, see e.g. Section 2.3.5
in \cite{Fe}, for every measurable function $\varrho:{\Bbb
R}^n\to[0,\infty]$, there is a Borel function $\varrho^*:{\Bbb
R}^n\to[0,\infty]$ such that $\varrho^*=\varrho$ a.e. in ${\Bbb
R}^n$. Thus, by Theorem \ref{th8.2.12}, $\varrho$ is measurable on
$p$-a.e. $k$-dimensional surface $S$ in ${\Bbb R}^n$ for every
$p\in(0,\infty)$ and $k=1,\ldots,n-1$. \end{rem}

\bigskip

We say that a Lebesgue measurable function $\varrho:{\Bbb
R}^n\to[0,\infty]$ is {\bf $p$-extensively admissible} for a family
$\Gamma$ of $k$-dimensional surfaces $S$ in ${\Bbb R}^n$, abbr.
$\varrho\in{\rm ext}_p\,{\rm adm}\,\Gamma$, if
\begin{equation}\label{eq8.2.15}\int\limits_S\varrho^k\,d{\cal A}\geqslant1\end{equation}
for $p$-a.e. $S\in\Gamma$. The {\bf $p$-extensive modulus}
$\overline M_p(\Gamma)$ of $\Gamma$ is the quantity
\begin{equation}\label{eq8.2.16}\overline M_p(\Gamma)=
\inf\int\limits_{{\Bbb R}^n}\varrho^p(x)\,dm(x)\end{equation} where
the infimum is taken over all $\varrho\in{\rm ext}_p\,{\rm
adm}\,\Gamma$. In the case $p=n$, we use the notations $\overline
M(\Gamma)$ and $\varrho\in{\rm ext}\,{\rm adm}\,\Gamma$,
respectively. For every $p\in(0,\infty)$, $k=1,\ldots,n-1$, and
every family $\Gamma$ of $k$-dimensional surfaces in ${\Bbb R}^n$,
\begin{equation}\label{eq8.2.17} \overline M_p(\Gamma)=M_p(\Gamma)\,.\end{equation}

\section{Ring $Q$-homeomorphisms and their properties}
Recall some necessary notions. Let $E\,,F\,\subseteq\,\mathbb R^n$
be arbitrary domains. Denote by $\Delta(E,F,G)$ the family of all
curves $\gamma:[a,b]\,\rightarrow\,\mathbb R^n$, which join $E$ and
$F$ in $G$, i.e. $\gamma(a)\,\in\,E\,,\gamma(b) \,\in\,F$ and
$\gamma(t)\,\in\,G\,$ for $a\,<\,t\,<\,b\,.$ Set $d_0\,=\,{\rm
dist}\,(x_0\,,\partial G)$ and let
$Q:\,G\rightarrow\,[0\,,\infty]\,$ be a Lebesgue measurable
function. Denote
\begin{equation*}
A(x_0,r_1,r_2) = \{ x\,\in\,{\mathbb R}^n : r_1<|x-x_0|<r_2\}\ ,
\end{equation*}
and
\begin{equation}\label{2.11}
S_i\,=\,S(x_0,r_i) = \{ x\,\in\,{\Bbb R^n} : |x-x_0|=r_i\}\,\,,\ \ \
i=1,2.
\end{equation}

\medskip
We say that a homeomorphism $f:G\to \mathbb R^n$  is \textit{the
ring $Q$-homeomorphism with respect to $p$-module at the point
$x_0\,\in\,G,$} ($1<p\leq n$) if the inequality
\begin{equation}\label{2.12}
\mathcal M_{p}\,\left(\Delta\left(f(S_1),f(S_2),
f(G)\right)\right)\,\ \leq \int\limits_{A} Q(x)\cdot
\eta^{p}(|x-x_0|)\ d\,x
\end{equation}
is fulfilled for any ring $A=A( x_0,r_1,r_2),$\,\, $0<r_1<r_2< d_0$
and for every measurable function $\eta : (r_1,r_2)\to [0,\infty
]\,,$ satisfying
\begin{equation}\label{2.13}
\int\limits_{r_1}^{r_2}\eta(r)\ dr\ \geq\ 1\,.
\end{equation}
The homeomorphism  $f:G\to \mathbb R^n$ is the \textit{ring
$Q$-homeomorphism with respect to $p$-module in the domain $G$}, if
inequality $(\ref{2.12})$ holds for all points $x_0\,\in\,G\,.$ The
properties of the ring $Q$-homeomorphisms for $p=n$ are studied in
\cite{SS10}.

The ring $Q$-homeomorphisms are defined in fact locally and contain
as a proper subclass of $Q$-homeomorphisms (see \cite{MRSY09}). A
necessary and sufficient condition for homeomorphisms to be ring
$Q$-homeomorphisms with respect to $p$-module at a point given in
\cite{Sal11}, asserts:

\medskip
\begin{propo}\label{th5.1}
{\em Let $G$ be a bounded domain in $\mathbb R^n$, $n\ge 2$ and let
$Q:\,G\,\rightarrow\,[0,\,\infty]\,$ belong to $L^1_{\rm loc}$. A
homeomorphism $f:\,G\to \mathbb R^n$ is a ring $Q$-homeomorphism
with respect to $p$-module at $x_0\in G$ if and only if for any
$0<r_1<r_2< d_0= dist\,(x_0,\partial G)$,
\begin{equation*}
\mathcal M_p\left(\Delta\left(f(S_1),\,f(S_2),\,
f(G)\right)\right)\,\le\,\frac{\omega_{n-1}}{I^{p-1}}\,\,,
\end{equation*}
where $S_1$ and $S_2$ are the spheres defined in (\ref{2.11})
\begin{equation*}
I\ =\ I(x_0,r_1,r_2)\ =\ \int\limits_{r_1}^{r_2}\
\frac{dr}{r^{\frac{n-1}{p-1}} q_{x_0}^{\frac{1}{p-1}}(r)}\,,
\end{equation*}
and $q_{x_0}(r)\,$ is the mean value of $Q$ over $|x-x_0|\,=\,r\,.$
Note that the infimum in the right-hand side of $(\ref{2.12})$ over
all admissible $\eta$ satisfying (\ref{2.13}) is attained only for
the function
\begin{equation*}
\eta_0(r)=\frac{1}{Ir^{\frac{n-1}{p-1}}
q_{x_0}^{\frac{1}{p-1}}(r)}\,.
\end{equation*}}
\end{propo}

In this sections we establish the relationship
between the ring and lower $Q$-homeomorphisms with respect to
$p$-module.

\medskip
\begin{theo}\label{th7.1}
{\em Every lower $Q$-homeomorphism with respect to $p$-module
$f:G\rightarrow G^*$ at $x_0\in G,$ with $p>n-1$ and $Q\in
L^\frac{n-1}{p-n+1}_{\rm loc}$, is a ring
$\widetilde{Q}$-homeomorphism with respect to $\alpha$-module at
$x_0$ with $\widetilde{Q}=Q^{\frac{n-1}{p-n+1}}$ and
$\alpha=\frac{p}{p-n+1}$.}
\end{theo}

\bigskip

\cc
\section{Lower $Q$-homeomorphisms and Orlicz-Sobolev classes}\label{9}

Let $D$ and $D'$ be two bounded domains in $\mathbb R^n$, $n\ge 2$
and $x_0\in D$. Given a Lebesgue measurable function $Q:D\to
[0,\infty]$, a homeomorphism $f:D\to D'$ is called the
\textit{lower $Q$-homeomorphism with respect to $p$-modulus at} $x_0$
if
\begin{equation}\label{eq6.1}
\mathcal M_p(f(\Sigma_\varepsilon))\ge \inf\limits_{\rho\in{\rm
ext_{p}adm}\;\Sigma_\varepsilon}
\int\limits_{A_\varepsilon(x_0)}\frac{\rho^p(x)}{Q(x)}dm(x),
\end{equation}
where
\begin{equation*}
A_\varepsilon(x_0)=\{x\in\mathbb R^n:
\varepsilon<|x-x_0|<\varepsilon_0\},\quad
0<\varepsilon<\varepsilon_0,\quad 0<\varepsilon_0<\sup\limits_{x\in
G}|x-x_0|,
\end{equation*}
and $\Sigma_\varepsilon$ denotes the family of all  spheres
centered at $x_0$ of radii $r$, $\varepsilon<r<\varepsilon_0$,
located in $D$.

\medskip

\begin{theo}{}\label{thOS4.1} {\it Let $D$ and $D'$ be domains in ${\Bbb R}^n$,
$n\geqslant3$, and let $\varphi:[0,\infty)\to[0,\infty)$ be
increasing with the condition (\ref{eqOS3.0a}). Then each
homeomorphism $f:D\to D'$ of finite distortion in the class
$W^{1,\varphi}_{\rm loc}$ is a lower $Q$-homeomorphism at every
point $x_0\in D$ with $Q(x)=K_{p,f}(x)$.} \end{theo}

\medskip

{\it Proof.} Let $B$ be a (Borel) set of all points $x\in D$ where
$f$ has a total differential $f'(x)$ and $J_f(x)\ne 0$. Then,
applying Kirszbraun's theorem and uniqueness of approximate
differential, see, e.g., 2.10.43 and 3.1.2 in \cite{Fe}, we see
that $B$ is the union of a countable collection of Borel sets
$B_l$, $l=1,2,\ldots\,$, such that $f_l=f|_{B_l}$ are bi-Lipschitz
homeomorphisms, see, 3.2.2, 3.1.4 and 3.1.8 in \cite{Fe}. With no
loss of generality, we may assume that the $B_l$ are mutually
disjoint. Denote also by $B_*$ the set of all points $x\in D$
where $f$ has the total differential but with $f'(x)=0$.

By the construction the set $B_0:=D\setminus \left(B\bigcup
B_*\right)$ has Lebesgue measure zero, see Theorem \ref{thOS4.0}.
Hence by Theorem 2.4 in \cite{KR$_2$} or by Theorem 9.1 in
\cite{MRSY} the area ${\cal A}_{S}(B_0)=0$ for a.e. hypersurface $S$
in ${\Bbb R}^n$ and, in particular, for a.e. sphere $S_r:=S(x_0,r)$
centered at a prescribed point $x_0\in\overline{D}$. Thus, by
Corollary \ref{corOS3.2} ${\cal A}_{S_r^*}(f(B_0))=0$ as well as
${\cal A}_{S_r^*}(f(B_*))=0$ for a.e. $S_r$ where $S_r^*=f(S_r)$.

Let $\Gamma$ be the family of all intersections of the spheres
$S_r$, $r\in(\varepsilon,\varepsilon_0)$,
$\varepsilon_0<d_0=\sup\limits_{x\in D}\,|x-x_0|,$ with the domain
$D.$ Given $\varrho_*\in{\mathrm adm}\,f(\Gamma)$,
$\varrho_*\equiv0$ outside of $f(D)$, set $\varrho\equiv 0$
outside of $D$ and on $B_0$,
$$\varrho(x)\ \colon=\ \varrho_*(f(x))\Vert f'(x)\Vert \qquad{\rm for}\ x\in D\setminus B_0\,.$$

Arguing piecewise on $B_l$, $l=1,2,\ldots$, we have by 1.7.6 and
3.2.2 in \cite{Fe} that $$\int\limits_{S_r}\varrho^{n-1}\,d{\cal A}\ \geqslant\
\int\limits_{S_{*}^r}\varrho_{*}^{n-1}\,d{\cal A}\ \geqslant\ 1$$
for a.e. $S_r$ and, thus, $\varrho\in{\mathrm{ext_p\,adm}}\,\Gamma$.

The change of variables on each $B_l$, $l=1,2,\ldots\,,$ see,
e.g., Theorem 3.2.5 in \cite{Fe}, and countable additivity of
integrals give the estimate
$$\int\limits_{D}\frac{\varrho^p(x)}{K_{p,f}(x)}\,dm(x)\ \leqslant\
\int\limits_{f(D)}\varrho^p_*(x)\, dm(x)$$ and the proof is
complete.

\medskip

\begin{corol}{}\label{corOS4.1} {\it Each homeomorphism $f$ of finite distortion
in ${\Bbb R}^n$, $n\geqslant3$, in the class $W^{1,s}_{\rm loc}$
for $s>n-1$ is a lower $Q$-homeomorphism at every point
$x_0\in D$ with $Q(x)=K_{p,f}(x)$.} \end{corol}

\medskip

\cc

\medskip

\end{document}